\documentclass[10pt,letterpaper]{article}
\usepackage[top=0.85in,left=2.75in,footskip=0.75in]{geometry}

\usepackage{amsmath,amssymb}

\usepackage{changepage}

\usepackage[utf8x]{inputenc}

\usepackage{textcomp,marvosym}

\usepackage{fixltx2e}

\usepackage{amsmath,amssymb}

\usepackage{cite}

\usepackage{nameref,hyperref}

\usepackage{microtype}
\DisableLigatures[f]{encoding = *, family = * }

\usepackage[table]{xcolor}

\usepackage{array}

\newcolumntype{+}{!{\vrule width 2pt}}

\newlength\savedwidth

\usepackage{rotating}

\usepackage{multirow}
\usepackage{booktabs}
\usepackage{tabularx}

\raggedright
\usepackage[aboveskip=1pt,labelfont=bf,labelsep=period,justification=raggedright,singlelinecheck=off]{caption}

\makeatletter
\renewcommand{\@biblabel}[1]{\quad#1.}
\makeatother

\date{}

\usepackage{lastpage,fancyhdr,graphicx}
\usepackage{epstopdf}
\fancyheadoffset[L]{2.25in}
\fancyfootoffset[L]{2.25in}
\usepackage{listings}
\begin{document}
\vspace*{0.2in}

\begin{flushleft}
{\Large
\textbf\newline{The Effect of Gender in the Publication Patterns in Mathematics}
}
\newline
\\

Helena Mihaljević-Brandt\textsuperscript{1,\P, *},
Lucía Santamaría\textsuperscript{2,\P, *},
Marco Tullney\textsuperscript{3,\P, *}

\bigskip
\textbf{1} Independent Researcher, Berlin, Germany
\\
\textbf{2} Independent Researcher, Potsdam, Germany
\\
\textbf{3} German National Library of Science and Technology (TIB), Hannover, Germany
\\
\bigskip

%
%
\P These authors contributed equally to this work.





* helena.m-b@t-online.de (HMB) \url{http://orcid.org/0000-0003-0782-5382}; lucia.santamaria@ymail.com (LS) \url{http://orcid.org/0000-0002-5986-0449}; marco.tullney@tib.eu (MT) \url{http://orcid.org/0000-0002-5111-2788}

\end{flushleft}
\section*{Abstract}
Despite the increasing number of women graduating in mathematics, a systemic gender imbalance persists and is signified by a pronounced gender gap in the distribution of active researchers and professors. Especially at the level of university faculty, women mathematicians continue being drastically underrepresented, decades after the first affirmative action measures have been put into place. A solid publication record is of paramount importance for securing permanent positions. Thus, the question arises whether the publication patterns of men and women mathematicians differ in a significant way. Making use of the zbMATH database, one of the most comprehensive metadata sources on mathematical publications, we analyze the scholarly output of $\sim\!150,000$ mathematicians from the past four decades whose gender we algorithmically inferred. We focus on development over time, collaboration through coautorships, presumed journal quality and distribution of research topics---factors known to have a strong impact on job perspectives. We report significant differences between genders which may put women at a disadvantage when pursuing an academic career in mathematics.

\section*{Introduction}

The overall number of women in Science, Technology, Engineering and Mathematics (STEM) has undoubtedly been increasing in the past decades~\cite{West.JD.Jacquet.J.ea:2013}. There is a growing presence of women on all academic levels: more female students, university graduates, PhDs, and professors. However, gender inequities persist, in particular towards higher academic ranks. Women are still significantly underrepresented when it comes to permanent academic positions such as tenure-track and university faculty posts.

Providing comprehensive figures and statistics on women in mathematics is not trivial, due to scattered data among countries, national differences in academic careers, and the difficulty to exactly define what ``mathematical research'' entails, among others. Among the available surveys,~Hobbs and Koomen~\cite{Hobbs.C.Koomen.E:2006} is the most complete report on women in mathematical research in Europe and shows a comparison between data taken in 1993 and 2005. While the proportion of female mathematicians had increased in the 12-year-long span in almost all countries, women are still underrepresented, especially in northern and western Europe. More specifically, the London Mathematical Society's report on Women in Mathematics from 2013 surveyed the university mathematics departments in the UK and states a 29\% female representation at the lecturer and a mere 6\% at the professor level~\cite{McWhinnie.S.Fox.C:2013}. For the faculty positions in the U.S., the Conference Board of Mathematical Sciences (CBMS) 2015 study~\cite{Blair.R.Kirkman.EE.ea:2013} reports an increase in percentages of women among tenured and tenure-eligible faculty, with a relatively low level of 11\% women in tenured faculty of doctoral university departments. On the other hand, the share of doctorates granted to women remained almost constant at around 30\% between 2000 and 2010~\cite{Blair.R.Kirkman.EE.ea:2013}. For a list of more descriptive statistics see~\cite{Frietsch.R.Haller.I.ea:2009}.

Scholarly publications are (still) the major means of communicating research results. A solid publication record is a key factor in achieving and sustaining a successful academic career in every discipline. The relentless pressure to promptly and steadily publish in the hope of sustaining one's career has lead to the coining of the phrase ``publish or perish''. The central role of scholarly productivity when pursuing a professorial post bears the question whether men and women in mathematics differ regarding their publication patterns, and whether these differences have changed over time. 

Due to the difficulty in obtaining complete records, many studies on gender effects in academia are restricted to small amounts of data, mostly from surveys, or limited to a few journals or countries. While small-data studies might suffer from lack of generality and the danger of ``anecdata'', they do typically deal with cleaner, less ambiguous data sources. In particular, surveys based on self-reported data have built-in authorship disambiguation and can acquire gender assignments from the respondents. Large-scale publication databases instead often suffer from missing information, which then needs to be inferred. In this article, we analyze publication metadata from zbMATH, an online bibliographic database comprising more than 3 million mathematical publications classified by research topics and enhanced with author profiles. Based on the authors' first names, we attribute a gender to $\sim\!150,000$ mathematicians from the past four decades. To the best of our knowledge, this is the most comprehensive large-scale study of this kind in the field of mathematics. 

Our analysis addresses four dimensions which are crucial for understanding the impact of publication patterns on the gender gap in mathematics: (1) research activity over time, (2) choice of journals and their perceived quality, (3) collaboration through coauthorship, and (4) distribution across specialisation fields. 

The growing presence of women in mathematics in particular and in STEM disciplines in general has often led to the belief that achieving a balanced gender ratio is merely a matter of time. Yet, the increasing number of women mathematicians entering graduate school will only translate into an overall increase of women in academia if they move on towards higher scholarly ranks, which motivates our longitudinal analysis. For active researchers in mathematics, journal articles are the predominant publication form, and the perceived quality of a journal is an important criterion when judging the excellence of a researcher's publication record. Said quality is often measured by means of journal rankings, which are highly controversial. Nevertheless, they are widely in use, often relieving evaluators of the onus of reading and judging an article by its content. For authors pondering which journal to submit to, or for university administrators, research funders, or governmental bodies evaluating publication lists, journal rankings are presumed to be helpful, thus we employ them in our study. A further career determinant worth analyzing is the interaction with other researchers: a prospering network advances the chance of citations by peers, facilites building collaborations and can help securing a university position. On the other hand, the share of single-authored publications in mathematics is, in contrast to many experimental sciences, very high. Such publications act as proof of one's individual capabilities, and are thus particularly valuable at an early career stage. They show that one is not ``completely dependent on senior people for ideas, guidance, techniques, […] Hence, one is ready for a faculty position''~\cite{McKensie:2012}. Last but not least, mathematical topics differ in terms of general popularity, preferences of journal editors, the overall community, and particular gatekeepers. Some topics bear higher chances for a scientific career in mathematics than others, and popular journals reflect a certain bias for particular mathematical fields~\cite{Mihaljevic-Brandt.H.Teschke.O:2014}.

Our analysis shows a higher dropout rate of women as authors of mathematics papers, despite an upwards trend in the share of women mathematicians over the last four decades. Women publish less articles than men during their careers and write less articles alone but establish networks of comparable size. Their work is significantly less published in journals of high presumed quality, measured in terms of the Thomson Reuter's journal impact factor (JIF) as the most prominent journal ranking schema and the manually compiled Australian ERA ranking. In three highly regarded journals, Annals of Mathematics, Inventiones Mathematicae, and Journal of the AMS, a significant underrepresentation of women persists without notable uptrend over the last four decades. While the presence of women is shown to strongly fluctuate across mathematical fields, the choice of research fields is found not to be responsible for the substantial underrepresentation of women in said three journals.

\subsection*{Related work}

Several studies have looked at publication patterns in the context of gender.  
Frietsch~et al.~\cite{Frietsch.R.Haller.I.ea:2009} studied women's share among patent applications and scientific publications. Using two large data sets from the 1990s and 2000s, the authors estimate the average proportion of women in mathematics publications of 14 countries in 2005 as 16.5\% in total, being lowest in Germany (11.7\%) and highest in Italy (19.6\%). Bentley~\cite{Bentley.P:2012} studied data from two international surveys in 1993 and 2007 involving academic staff in Australian public universities. While gender differences in publication productivity decreased significantly within the 14-year time period, in 2007 women published on average only 76\% of what men published.

West~et al.~\cite{West.JD.Jacquet.J.ea:2013} analyzed JSTOR data, a digital library corpus spanning several centuries of publications in sciences and humanities. Across all disciplines, the authors found 21.9\% women; for mathematics, the estimated percentage of female authors is significantly lower at only 10.6\% (= 6,134 persons). The authors showed that authorships by women are not evenly distributed over the dimensions time, field/subfield, coauthorship, and authorship position, the latter being particularly relevant for disciplines such as molecular biology or medicine. A recent, large-scale study based on Thomson Reuters' Web of Science (WoS) unveiled the persistence of women's underrepresentation in Russian science between 1973 and 2012 across many disciplines, particularly in mathematics and physics~\cite{Paul-Hus.A.Bouvier.RL.ea:2015}. Using fractional authorships computed by weighting an authorship as a reciprocal of the total number of (gender-identified) authors of the respective publication, the authors show disparities between women and men regarding research output, productivity, collaboration and scientific impact. Employing similar methods, the analysis of all articles published between 2008 and 2012 and indexed in WoS by Larivi\`{e}re~et al.~\cite{Lariviere2013} confirms the persistence of gender imbalances in research output worldwide.

Collaboration, e.g. in the form of joint research, publishing or grant proposals, is considered a strong indicator for productivity and integration within the research community, and hence has been the subject of various studies. For various career stages of mathematicians, a positive correlation between the number of distinct coauthors and the overall number of publications has been observed~\cite{Hu.Z.Chen.C.ea:2014}. 
The understanding of the structural factors that affect men's and women's collaboration could lead to the development of policies to address differences in power and resource dynamics~\cite{Bozeman.B.Gaughan.M:2011}. A similar argument is made regarding interdisciplinary research and the possibility to use it to foster career success. Bosquet and Combes~\cite{Bosquet.C.Combes.P:2013} use publication records of a group of French economists to evaluate the impact of various factors on two measures of article quality defined via journal quality and Google Scholar citation numbers, respectively. Among others, they show a high impact of the coauthor network size with respect to both quality measures. Effects of topical diversity and specialisation are considered as well, showing disparities between single subfields regarding the number and average quality of publications, as well as the number of citations.

According to Garg and Kumar~\cite{Garg.KC.Kumar.S:2014}, female and male Indian researchers in life sciences differ with respect to team work and choice of journals: women tend to work in smaller teams with less international collaborative papers, and their publications are more often published in domestic journals with lower impact factors. Using a set of high-quality Spanish journals across various fields of science,~Maule{\'o}n~et al.~\cite{Mauleon.E.Hillan.L.ea:2013} analyzed the gender distribution over time of authors and editorial board members as representatives of active and leading scientists, respectively. While the share of female authors grew over the years, only a minimal effect was visible in editorial board membership. Restricting to eight international mathematics journals,~Maule{\'o}n and Bordons~\cite{Mauleon.E.Bordons.M:2012} calculated the share of women in 2004 to be 4\% in editorial boards and 10\% of all authors. The authors could not discover any relation between gender and the journals' impact factors. They showed only slight differences in the collaborative behaviour and a slight disadvantage for women regarding citation metrics. 

Finally, questions related to productivity, amount of scholarly output, and merit are the topic of multiple studies. For a discussion on why some faculty members in political science publish more than others, see Hesli and Lee~\cite{Hesli.VL.Lee.JM:2011}. The construction of the concept of ``academic excellence'' is analyzed in Brink and Benschop~\cite{Brink.M.Benschop.Y:2012}. Recently, in a detailed discussion on productivity based on arXiv publications, Pierson~\cite{Pierson:2014} expresses hope that the gender productivity gap is closing, albeit recognizing that there is still a long way to go until gender parity is accomplished. Her analysis of almost 1 million arXiv articles on STEM disciplines supports some of the findings that we present in this paper.

\section*{Materials and methods}

\subsection*{The zbMATH corpus}

zbMATH (formerly Zentralblatt MATH) is an abstracting and reviewing service in mathematics, produced by FIZ Karlsruhe – Leibniz Institute for Information Infrastructure and coedited by the European Mathematical Society (EMS) and the Heidelberg Academy of Sciences and Humanities. Along with MathSciNet, zbMATH is one of the main services that provide bibliographic information (e.g. author names, title, source), abstracts and expert reviews as well as aggregated data (e.g. keywords) on mathematical publications to the scientific community. 

\paragraph*{Author profiles}\label{subsubsec:author_profiles}

The process of attributing authorship of bibliographic records to specific individual researchers, known as ``author name disambiguation'', is essential to perform analyses on publication data. This is by no means a straightforward task, due to the fact that personal names are not sufficiently distinct to unambiguously identify each individual within the large number of researchers active today and in the past. The problem is exacerbated by the inconsistent way in which author names appear in publications.

Author name disambiguation, and the more general problem of word sense disambiguation, is an open issue in computational linguistics and natural language processing. A taxonomy of automatic author name disambiguation methods can be found in~\cite{Ferreira:2012:BSA:2350036.2350040}. zbMATH provides author profiles for its indexed publications. They are constructed following a multifaceted approach based on automatic, manual, and collaborative methods. Initially, the assignment of a publication to an author is performed algorithmically, using mainly the name string and coauthorship information; some authorship assignments are post-processed manually by the zbMATH staff, and, recently, users can contribute corrections through an online form~\cite{DBLP:conf.mkm.Mihaljevic-BrandtMR14}.

Name-string-based authorship attribution suffers from various well-known issues such as confusion of authors with frequent names (John Miller), missing name parts (middle names), variability through transliteration (Chebyshev vs. Tschebyscheff), or name changes throughout a researcher's life (after a change of civil state). In zbMATH, an ambiguous authorship is initially attributed to all possible candidates, e.g. an article authored by ``P. Smith'' will appear in the profiles of ``Paul Smith'' and ``Patricia Smith''. Hence an author profile in zbMATH can contain wrongly assigned articles, or publications can be missing. A combination of both issues is possible as well. Some types of assignment error are more frequent in relation to certain regions; for instance, authors with Chinese names are often assigned to ambiguous profiles due to common surnames and transliteration problems. Due to this latter factor, authors from the former Soviet countries frequently suffer from profile splitting.

Cases such as ``P. Smith'' referred above, when the author name string does not suffice to unequivocally assign an author, are marked as ambiguous and are not to be trusted. We excluded such authorships from our analyses.

\paragraph*{Mathematics Subject Classification (MSC)}

The Mathematics Subject Classification (MSC) is a universally accepted classification schema in mathematics jointly produced and maintained by MathSciNet and zbMATH. It is a tree-like, three-level alphanumerical scheme to label publications according to their subject matter. For instance, ``35'' stands for the area of ``Partial Differential Equations'', ``35J'' refers to ``Elliptic Equations and Systems'', and ``35J75'' is the classification for ``Singular Elliptic Equations''. The MSC exists since 1970; it is revised every 10 years, with MSC 2010 being the current version~\cite{MSC:2010}. All entries in zbMATH since 1970, and even some older data, are classified by its editorial board according to the MSC. The first code of an article is its primary classification, supposedly reflecting its most important content; subsequent codes are secondary.

Articles classified with first-level MSC codes between 03 and 65 belong to mathematics subfields such as Logic, Discrete Mathematics, Algebra, Analysis, Geometry, Topology, Probability Theory, Statistics, and Numerical Mathematics. Those with codes between 68 and 94 fall into mathematical applications related to Computer Science, Physics, or Economics. The areas 00 and 01 refer to general topics and history of mathematics; 97 is Mathematical Education.

\paragraph*{Focus on core mathematicians}\label{subsubsec:core_mathematicians}

The zbMATH corpus includes publications in which mathematical concepts are applied to scientific subjects such as quantum theory or computational biology. Hence, the zbMATH author records necessarily contain individuals whose primary research is not in mathematics. Their publications usually appear in interdisciplinary journals. For our analysis we have decided to only consider authors as mathematicians if they have published at least one article in a ``core mathematics journal''. To extract such journals, we resort to the MSC labels applied to a journal's articles to identify core mathematics publications, using the following heuristics: a journal belongs in the group of core mathematics journals if at least 90\% of its articles with primary classification between 03 and 94 fall within the range 03–65. Our corpus consists of publications drawn from 4064 journals. By including the rather applied areas of Statistics (62) and Numerical Mathematics (65) we retain a comparatively large number of 1,050 core mathematics journals, representing 26\% of all serials indexed by zbMATH with publications after 1970.

\subsection*{Guessing an author's gender}\label{author_gender}

In general, the gender of a person cannot be inferred programmatically from the name string---hence the use of a database of annotated names for gender assignment is often unavoidable. Several studies have made use of annotated data from varied sources:~West~et al.~\cite{West.JD.Jacquet.J.ea:2013} used a two-gender data set provided by the US Social Security Administration containing the top 1,000 annual girl and boy baby names since 1888; Frietsch~et al.~\cite{Frietsch.R.Haller.I.ea:2009} worked with a first name database comprising 8,291 names from six European countries. This database was created as part of a project carried out by the European commission in 2002, in which first name assignments were collected from various data sources, then post-processed in order to achieve a high quality level, and finally tested~\cite{Naldi.F.Luzi.D.ea:2002}. Implementations of gender guessing procedures include the Python library \textit{Gender Detector}~\cite{Vanetta.M:github:2014} and the R package \textit{gender}~\cite{Mullen.L:2016}. A benchmark comparison of gender detection algorithms can be found at \cite{Vanetta.M:2014}. 

We tried to attribute a gender to all zbMATH author profiles. The zbMATH corpus contains no supplementary metadata on authors' genders. We resorted to the annotated collection of first names described in~\cite{Michael:2007}. The dictionary file contains a list of more than 42,000 names with gender assignment~\cite{Michael:2008}; at the time of the last update in November 2008 the list was assumed to cover the vast majority of names in all European and some non-European countries like China, India, Japan, USA. According to the author Jörg Michael, for some countries such as Turkey or Korea all names have been independently classified by native speakers, indicating a high quality of the data. Yet the zbMATH data contains no information on the country of origin of the authors, hence we could not use this additional facet. Every name in the dictionary is labelled as ``unknown'', ``female'' or ``male'', with additional information on the certainty levels of the assignment.

On top of this algorithmic assignment, we have performed a manual check of the most prolific profiles. Specifically, we have reviewed the gender assignment for the top (in number of publications) 300 profiles classified as men and 500 women. These typically correspond to well-known mathematicians, whose gender we were able to assess. We found no instances of women mathematicians being incorrectly classified and a $\sim\!5\%$ error in men that were incorrectly labelled as women. We manually changed the gender of the latter, resulting on cleaner and more trustworthy gender data. We observed that these errors were more frequent for names from India, Pakistan and Japan. We conclude that our study overrepresents the amount of women by at most $\sim\!5\%$.

We point out that troughout the text, we speak of men and women as if those categories were real. Gender binary is a concept that lacks both complexity and suitability. But categories are constructs and results of societal action. Almost all databases and publications connecting names with gender operate within this binary schema. Since we do not have access to detailed author profiles and data on how authors identify when it comes to gender, we use the name as an indicator of the perceived gender of each author. Thus, our categorization of mathematicians into men and women is following societal and academic action and is due to pragmatic and substantive reasons.

\subsection*{Data overview}\label{subsec_statistical_description}

Our data set contains information on all publications by authors we identified as ``core mathematicians'' as described in subsection \textit{Focus on core mathematicians}. For every such publication we extracted the authors, their order of appearance, the title, the publication year, all MSC classification codes, and the journal name together with ISSN and e-ISSN codes, if the publication appeared in a journal. We additionally tagged authorships with a value indicating the quality of the authorship assignment (manually curated, algorithmically unambiguous, and ambiguous). The MSC classification exists since 1970. We restricted our analyses to authorships published since then to ensure consistency when addressing our key questions.

The corpus consists of 2,245,205 document entities corresponding to 3,433,553 instances of authorship, yielding 1.53 authors per article. The amount of ambiguous authorships which could not be assigned by zbMATH to a unique author amounts to 226,500 ($6.6\%$). We excluded these authorships from the analyses since we cannot decide which attribution might be correct.

We applied a gender assignment algorithm as described in~\cite{mihaljevic_brandt_2016_51325} to 262,146 author identities, the remaining 47,319 ($\sim\!15\%$) authors were listed without a first name. We were able to assign a gender to $\sim\!55\%$ of all profiles with a real first name; among those, 27,596 ($\sim\!19\%$) authors were classified as women and 116,657 as men.

\begin{figure}[h!]
	\begin{center}
		\includegraphics[width=\columnwidth]{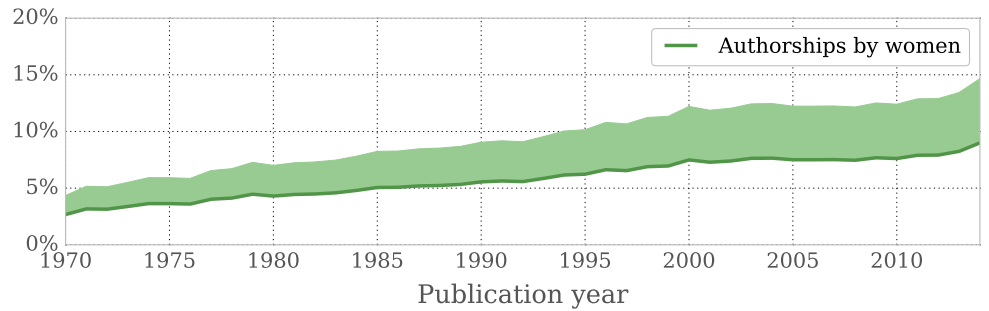}
		\caption{{\bf Percentage of authorships by women since 1970.} 
		The darker green line shows the share of authors that could be  identified as women, marking the lower limit. The shadowed region delimits the extrapolated values taking into account that we could only assign a gender to $61.4\%$ of all authorships. We estimate the proportion of women in the remaining group of authors to be smaller than their share in the identified group and conclude that the percentage of women among all authorships lies within the shadowed region.}
		\label{Fig 1}
	\end{center}
\end{figure}

We were able to guess the gender for $61.4\%$ of all (non-ambiguous) instances of authorship.  Fig~\ref{Fig 1} displays the percentage of all authorships that could be labelled as women since 1970. It is undeniable that the amount of women that publish in mathematics has been increasing---their share has tripled over the past four decades. 

As shown in Table \ref{Tab 1}, profiles of women are on average slightly less ambiguous than those of men; profiles with a high level of  ambiguity are significantly more frequent among those authors we could not classify as women or men. One explanation for the latter observation is the fact that our gender assignment algorithm requires the first name (which is missing  for $15.3\%$ of authors in our corpus), and the name string is also the main information used for authorship assignment by zbMATH. 

\begin{table}[h]
\begin{tabularx}{\textwidth}{ X X X X X X }
	\toprule
	 & \multicolumn{5}{ c }{Ratio of non-ambiguous publications per profile} \\
	 \cmidrule(r){2-6}
	 & 100\% & 90\% & 80\% & 70\% & 60\%\\
	 \midrule
	all authors & 49\% & 52\% & 55\% & 58\% & 65\% \\
	women & 71\% & 73\% & 75\% & 77\% & 84\% \\
	men & 62\% & 67\% & 71\% & 73\% & 81\% \\
	\bottomrule
\end{tabularx}
\caption{Distribution of the level of correctness of author profiles with respect to gender: to 71\% of all women and 62\% of all men no ambiguous authorships are assigned.}
\label{Tab 1}
\end{table}

\section*{Results}\label{results}

\subsection*{Cohort analysis}

A fundamental metric to understand the gender gap among scientists is the distribution of researchers that successfully sustain a career in academia over time. We analyse this aspect with a cohort analysis of authors since 1970. A cohort is a group of individuals that share a common characteristic over a defined period of time. Comparing the development of cohorts allows for longitudinal studies that unveil patterns over large periods of time. In our study, we assign authors to a given cohort according to the year of first recorded publication, which acts as an approximation for the beginning of their scientific career. We then follow the publication records of all authors in every cohort and select those individuals that continue publishing at least 5 and 10 years, respectively, after their first paper. We roughly assume researchers to be at the postdoctoral stage 5 years after their first publication and to have secured a permanent position at professorial level another 5 years later.

\paragraph*{Women mathematicians have tripled their number since 1970}

\begin{figure}[h]
	\begin{center}
		\includegraphics[width=\columnwidth]{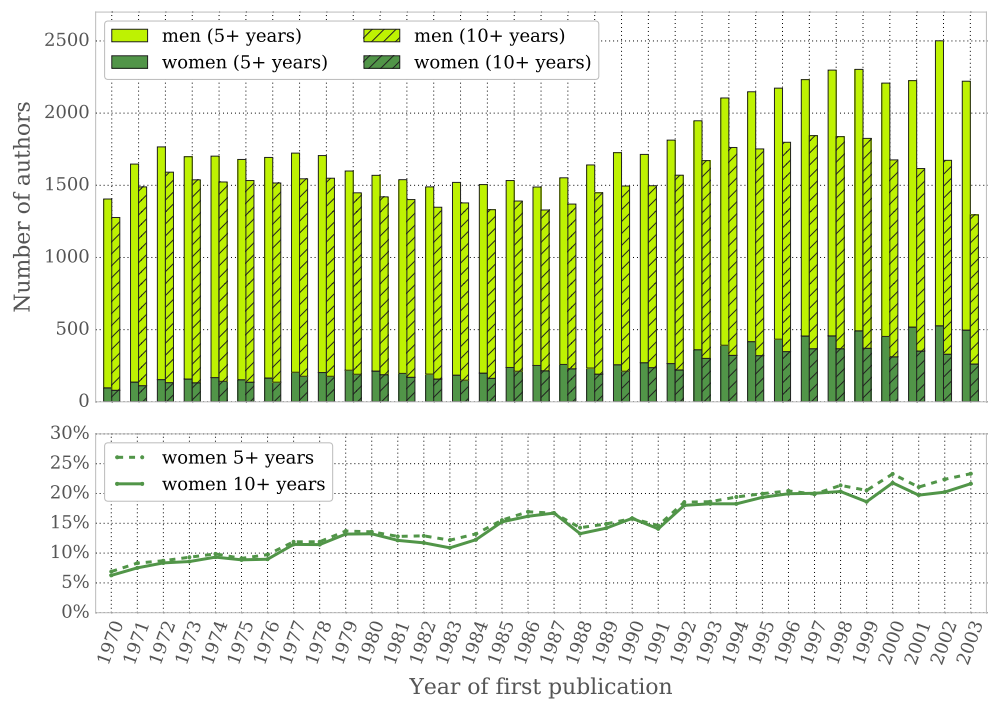}
		\caption{{\bf Cohort analysis since 1970.}
		Data corresponds to mathematicians with a first publication in the years 1970–2003 and another publication at least 5 and 10 years later, respectively. The upper panel shows the total number of authors per year and gender with continued careers 5 and 10 years after their first recorded publication. The lower panel displays the percentage of authors identified as women per cohort year and with careers lasting at least 5 and 10 years.}
		\label{Fig 2}
	\end{center}
\end{figure}

The total number of authors that continue publishing 5 and 10 years after their first paper appears in the upper panel of Fig~\ref{Fig 2}. Despite a brief stagnation in the 1980s, the number of mathematicians at our estimate of professorial level (10 years) has increased from $\sim\!1500$ in the 1970s to $\sim\!1800$ per year in the 1990s. Note that the decreasing number of authors after the 1990s reflects the fact that data for the most recent cohorts is partially incomplete: members of the early 2000s cohorts may well be still active today, eventhough they do not have a publication in the early 2010s. Therefore we regard the figures from the 2000s as non-conclusive.

Whereas a mere $\sim\!6\!-\!8\%$ of all authors in the early 1970s were women, their percentage has tripled within the following three decades (Fig~\ref{Fig 2}). Roughly 20\% of all authors with careers lasting at least 10 years that started publishing in the late 1990s are women. However, this is still far from the percentage of women that pursue mathematics education at the undergraduate level.

\paragraph*{Women publish less than men at the beginning of their careers}

\begin{figure}[h]
	\begin{center}
		\includegraphics[width=\columnwidth]{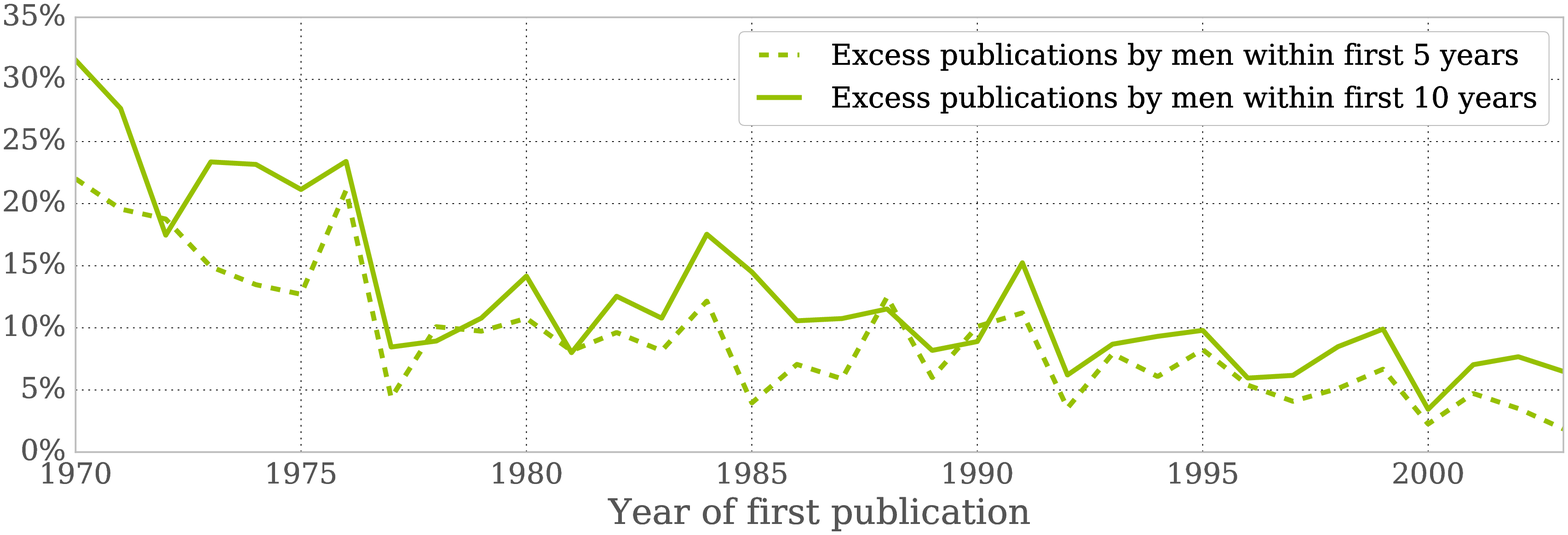}
		\caption{{\bf Excess of publications by men during the 5 and 10 years after their first article.} 
		On average, men publish 9\% and 13\% more papers than women within the first 5 and 10 years of their publishing career, respectively. Notably, the productivity gap has been closing since 1970.}
		\label{Fig 3}
	\end{center}
\end{figure}

An important metric that reveals differences between genders is the amount of publications. A shorter publication record at the time of applying for an academic position diminishes the chances of success. Our findings confirm the phenomenom known as productivity gap, agreeing with those of Pierson~\cite{Pierson:2014}. We have analyzed the average amount of publications of men and women during the 5 and 10 years after their first article appeared; we find that men routinely surpass women. Fig~\ref{Fig 3} illustrates this fact, showing the excess of publications by men per cohort year. Despite the visible trend towards a narrower productivity gap, men are shown to publish on average 9\% and 13\% more papers than women within their first 5 and 10 years as authors, respectively.

\paragraph*{Women leave academia at a higher rate than men}

\begin{figure}[h]
	\begin{center}
		\includegraphics[width=\columnwidth]{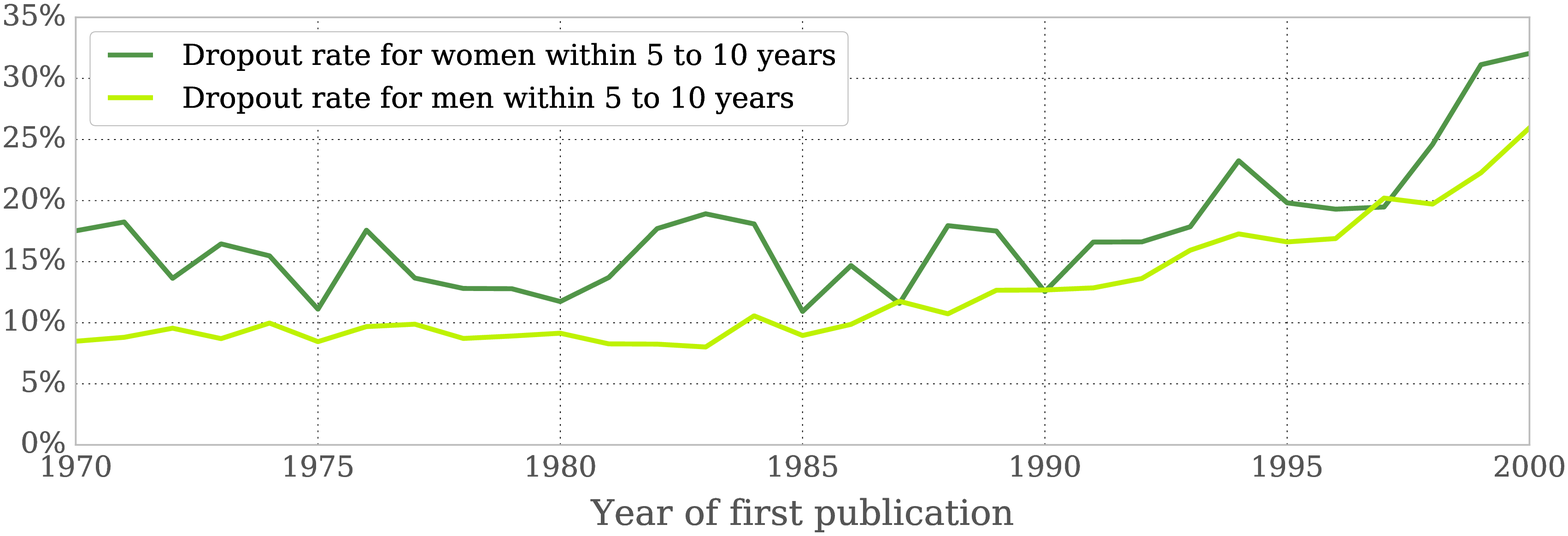}
		\caption{{\bf Dropout rate in authors per cohort year within 5 to 10 years of their first publication.}
		The figure shows the percentage of men and women that, having published for at least 5 years, disappear from our records 5 to 10 years after their first article.}
		\label{Fig 4}
	\end{center}
\end{figure}

In Fig~\ref{Fig 4} we evaluate the percentage of men and women that stop publishing at some point between 5 and 10 years after their first paper. This we interpret as them having discontinued their research careers. The results clearly show an overall increase in the dropout rate for both genders over the past decades. This phenomenon reflects the fact that an increasing number of doctorates and postdoctoral scholars find it impossible to secure a permanent position in science, thus being forced to leave academia even after a comparatively large publishing trajectory of at least 5 years.

On top of that, women end their research career in this period at a higher rate than men, as can be seen in Fig~\ref{Fig 4}. A larger comparative drop in the years that roughly correspond to the postdoctoral stage may indicate that women either lack support to reach the stable professorial ranks, or do not wish or are not able to continue holding temporary positions in a period of life often coincident with family-planning decisions. This evidences the existence of a leaky pipeline, understood as a continuous loss of women in academia as they climb the career ladder, that prevents the number of women in mathematics from increasing as fast as the rates of women enrolling at the graduate level would suggest.

\subsection*{Perception of quality}

We analyze two of the most prominent journal rankings in mathematics, the Australian ERA list and the journal impact factor. ERA (Excellence in Research for Australia) is an initiative of the Australian Research Council to measure the quality of Australian research. ERA compiled a ranked journal list with ranks reaching from A* (with the expectation that almost every article in these journals will be of very high quality) over A (most articles of very high quality) and B (only some articles of very high quality) to C (very high quality not expected) and ``not ranked''. We use the ERA 2010 ranked journal list which is the most recent list of this kind. It comprises 20,712 journals of all academic fields. Since 2012 ERA does not rank the journals anymore and publishes instead an unranked but more complex list of journals.

The journal impact factor (JIF) is calculated based on Thomson Reuter's Journal Citation Reports. It is defined as the average number of article citations in a journal in the last two years. There is widespread criticism of the impact factor, as there are many means to manipulate it~\cite{Kiesslich.T.Weineck.SB.ea:2016}. Based on different publication and citation practices in various academic fields, it is difficult to compare the impact factor across disciplines. We refer to the 2015 edition of Journal Citation Reports (which uses publication data from 2013/2014) and investigate impact factors for 7 subject categories: ``Mathematics'', ``Mathematics, Applied'', ``Mathematics, Interdisciplinary Applications'', ``Computer Science, Theory and Methods'', ``Logic'', ``Physics, Mathematical'', ``Statistics and Probability''.

We identify the ERA rank and JIF by ISSN, which leads to 1,768 matched journals in the ERA list and 728 journals with a JIF (663 journals appear in both indexes). We investigate publications from 2004–2013. Note that we refer to the most recent impact factor from 2015 which in most cases is not the impact factor of the journal at the submission or publication date. We further restrict our sample to publications with conclusive gender attribution and in journals from the ERA list or with an impact factor. We do not limit the analysis to core mathematics journals (see subsection \textit{Focus on core mathematicians}) after our tests showed no major differences with respect to the full analysis. 

For both rankings we investigate how the percentage of women among the authors per rank deviates from the average in the complete data set. We use the ranks A*, A, B, and C as assigned by ERA. Since the impact factor is not discrete, we calculate percentiles of the journals in our list, the highest (lowest) containing the 25\% journals with the highest (lowest) impact factors. We end up with four groups of publications for each ranking. In both cases, the number of articles per group differs. Furthermore, we examine authorships of single-authored publications versus all authorships (including single- and multi-authored) and quote separate results.

\paragraph*{High-ranked mathematics journals publish less articles authored by women}

\begin{figure}[h]
	\begin{center}
		\includegraphics[width=\columnwidth]{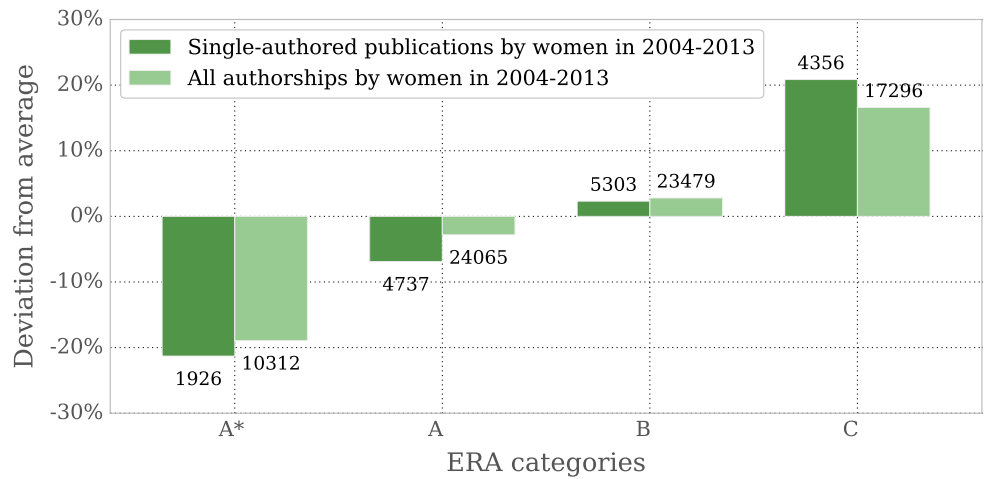}
		\caption{{\bf Deviation of the share of women's publications from the average ERA ranking.}
		ERA values rank from A* (excellent) to C (very high quality not expected). Annotated over or under each bar is the absolute number of single-authored publications and all authorships by female authors per group.}
		\label{Fig 5}
	\end{center}
\end{figure}

We uncover significant differences between female and male authors in both these rankings. There are fewer women amongst the authors of top-ranked journals than could have been expected from the overall distribution of men and women in our samples. Our analysis of the ERA ranked journal list shows that only $8.5\%$ women are authors of single-authored publications in A* journals, a deviation of $-21.3\%$ from the average, while there is a deviation of $+20.9\%$ in the low-ranked category C. The results for all single- and multi-authored papers are similar, since the A* journals show $18.9\%$ less authorships by women while there is an excess of $16.6\%$ with respect to the average in C journals. The results are displayed in Fig~\ref{Fig 5}.

\begin{figure}[h]
	\begin{center}
		\includegraphics[width=\columnwidth]{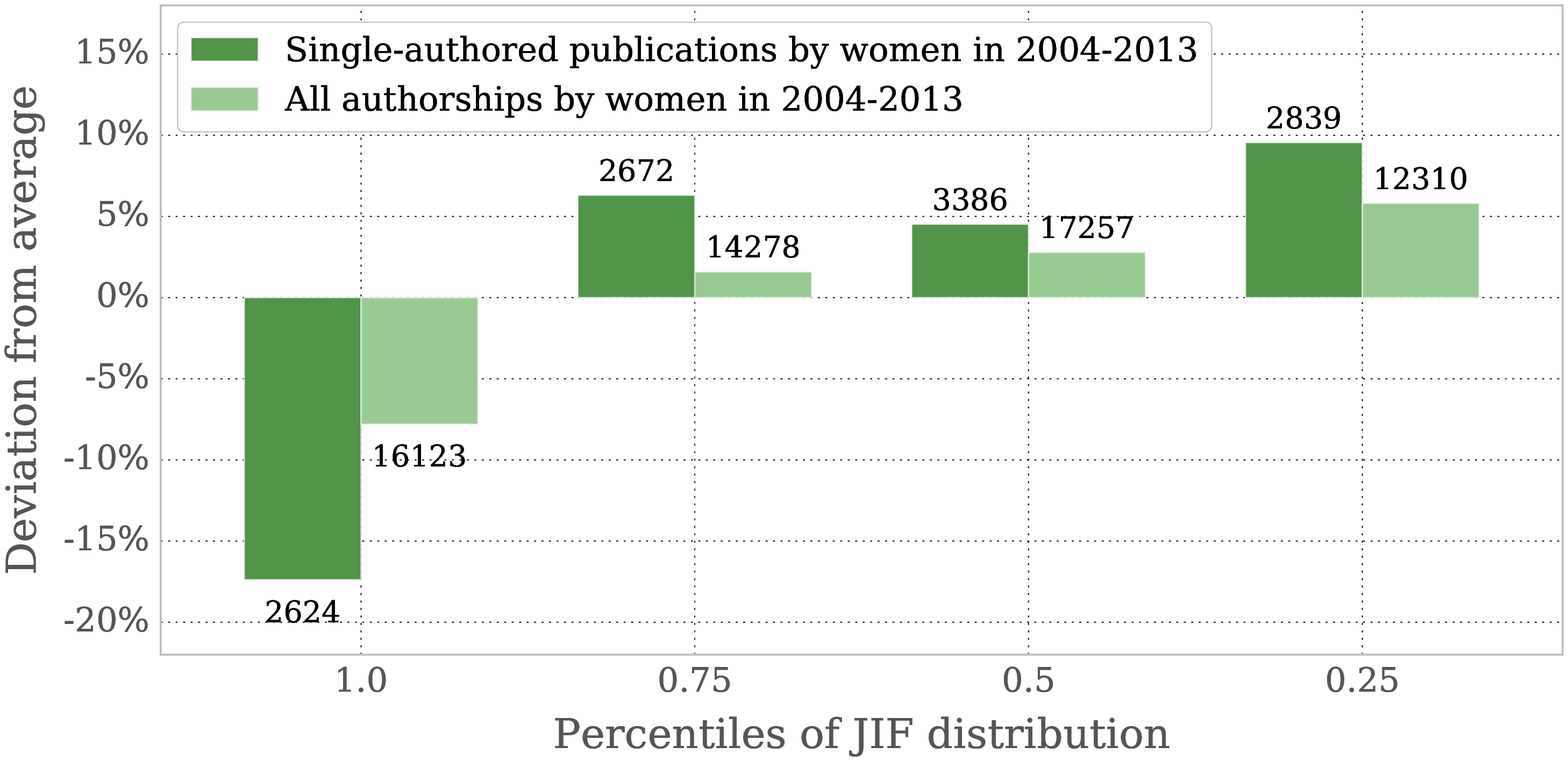}
		\caption{{\bf Deviation of the share of women's publications from the average impact factor.}
		JIF values are shown as percentile, each containing 25\% of the indexed journals in our sample. 1.0 to 0.75 contains the highest impact factors, 0.25 to 0.0 contains the lowest. Annotated over or under each bar is the absolute number of respective authorships by female authors per percentile bracket.}
		\label{Fig 6}
	\end{center}
\end{figure}

For the percentiles of the JIF we obtain a significantly smaller share of female authors in the top 25\% journals with $17.4\%$ and $7.8\%$ less women in the single-authored and all-authorships samples, respectively, when compared to the corresponding overall shares. Details are shown in Fig~\ref{Fig 6}.

\paragraph*{Women are drastically underrepresented in three top-ranked journals}

\begin{figure}[h!]
	\begin{center}
		\includegraphics[width=\columnwidth]{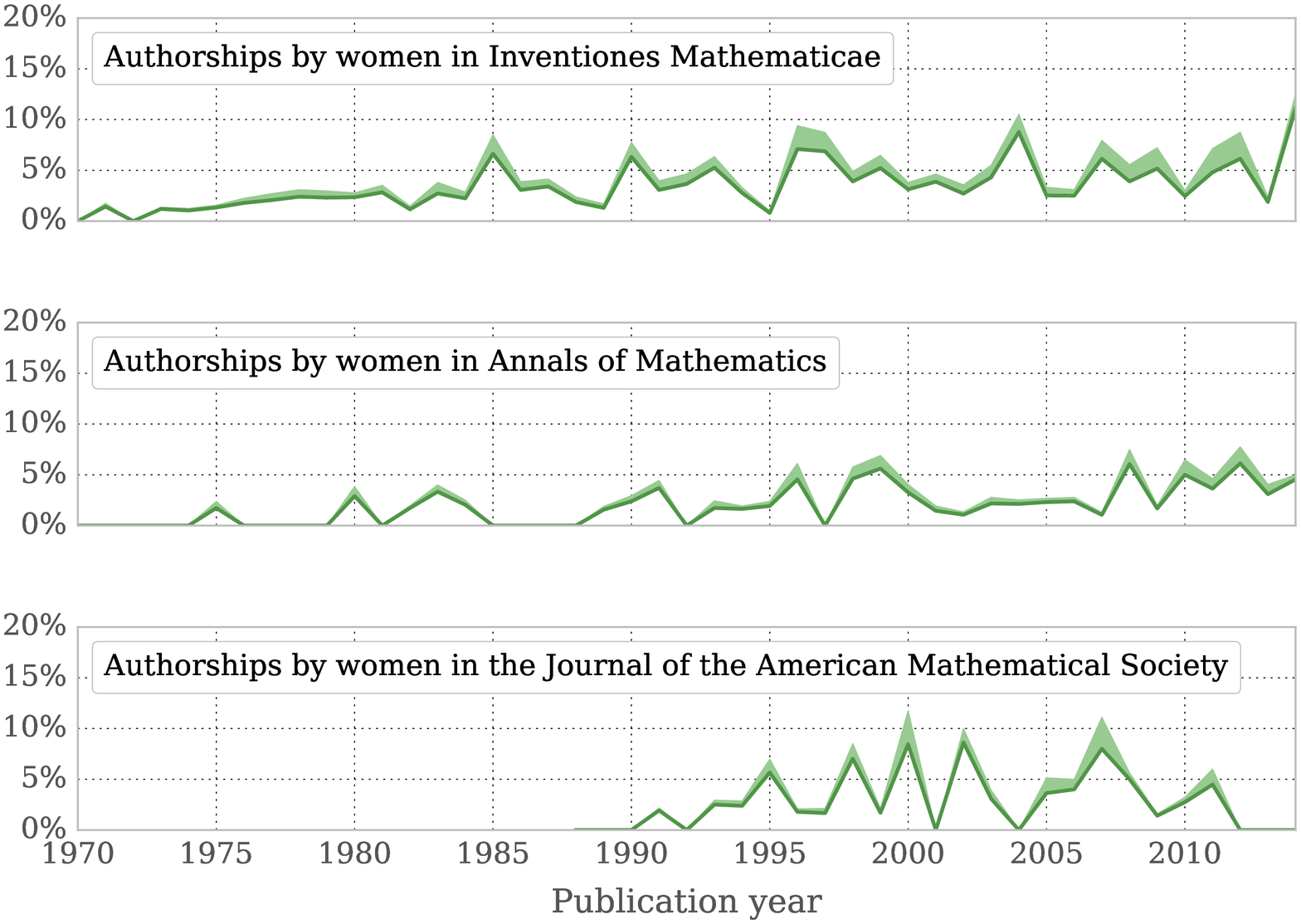}
		\caption{{\bf Share of authorships by women in three reference journals since 1970.} 
		Percentage of authorships by women in Inventiones Mathematicae, Annals of Mathematics, and Journal of the American Mathematical Society. The dark green line shows our estimation of the share of female authors, the shadowed region indicates the error margin we consider to be realistic (same proportion of women in the unknown group).}
		\label{Fig 7}
	\end{center}
\end{figure}

Female authors publish significantly less in the top journals according to the two ranking mechanisms we examined. To look further at the peak of what is perceived as good journal quality, we analyze the gender distribution of authors in three highly regarded serials: ``Annals of Mathematics'', ``Inventiones Mathematicae'', and ``Journal of the AMS''. Our aim is to shed light on the percentage of female authorships and its development over the last decades. Fig~\ref{Fig 7} illustrates that the share of female authors is significantly lower than in the overall sample of mathematical journals, but also lower than in the top categories of the ERA list and the JIF. Even more remarkable: no upward trend is discernible: in the Journal of the AMS, for instance, we do not identify a single female author in the last years (2012: 42 men/23 unknown; 2013: 43 men/16 unknown; 2014: 20 men/9 unknown). Note that these numbers include all authorships. Limiting the analysis to single-authored publications results in very low absolute numbers, but following a similar tendency.

\subsection*{Collaboration through coauthorship}

Often, the idea of mathematics is associated with the image of a solitary researcher, whose work requires only pen and paper and hence can be performed in pure isolation. It is therefore no surprise that more than 95\% of all publications of ``core mathematicians'' in zbMATH before 1950 were single-authored; in our corpus they still comprise almost 60\% of all records. Despite the strategical relevance of single-authored publications, a strong shift towards collaborative authoring is evident. This trend is supported by our data: among the publications that appeared after 2000, only 39\% had been authored by a single person. On the other hand, almost 25\% of the publications were authored by more than 2 authors compared to e.g. less than 8\% in the period 1970–2000, indicating a growing interest in building networks. 
We investigate coauthorship patterns in terms of two metrics: share of single-authored publications per author and number of distinct coauthors (``network size''). In contrast to other studies, we do not analyze the positions of coauthors with respect to gender, since in mathematics, authors are usually listed in alphabetical order (see, e.g.,~\cite{West.JD.Jacquet.J.ea:2013}).

We examine the records of mathematicians that started publishing in 1970 or later for whom we were able to successfully assign a gender (25,050 women and 77,792 men). Since the odds of collaborating with different coauthors increase with publication record length, we segment authors according to their number of publications. We sample the data as to include the same number of female and male authors per total number of publications. The sample size is thus determined by the maximum number of authors in the women's group. 

\paragraph*{Women publish less single-authored papers}

\begin{figure}[h!]
	\begin{center}
		\includegraphics[width=\columnwidth]{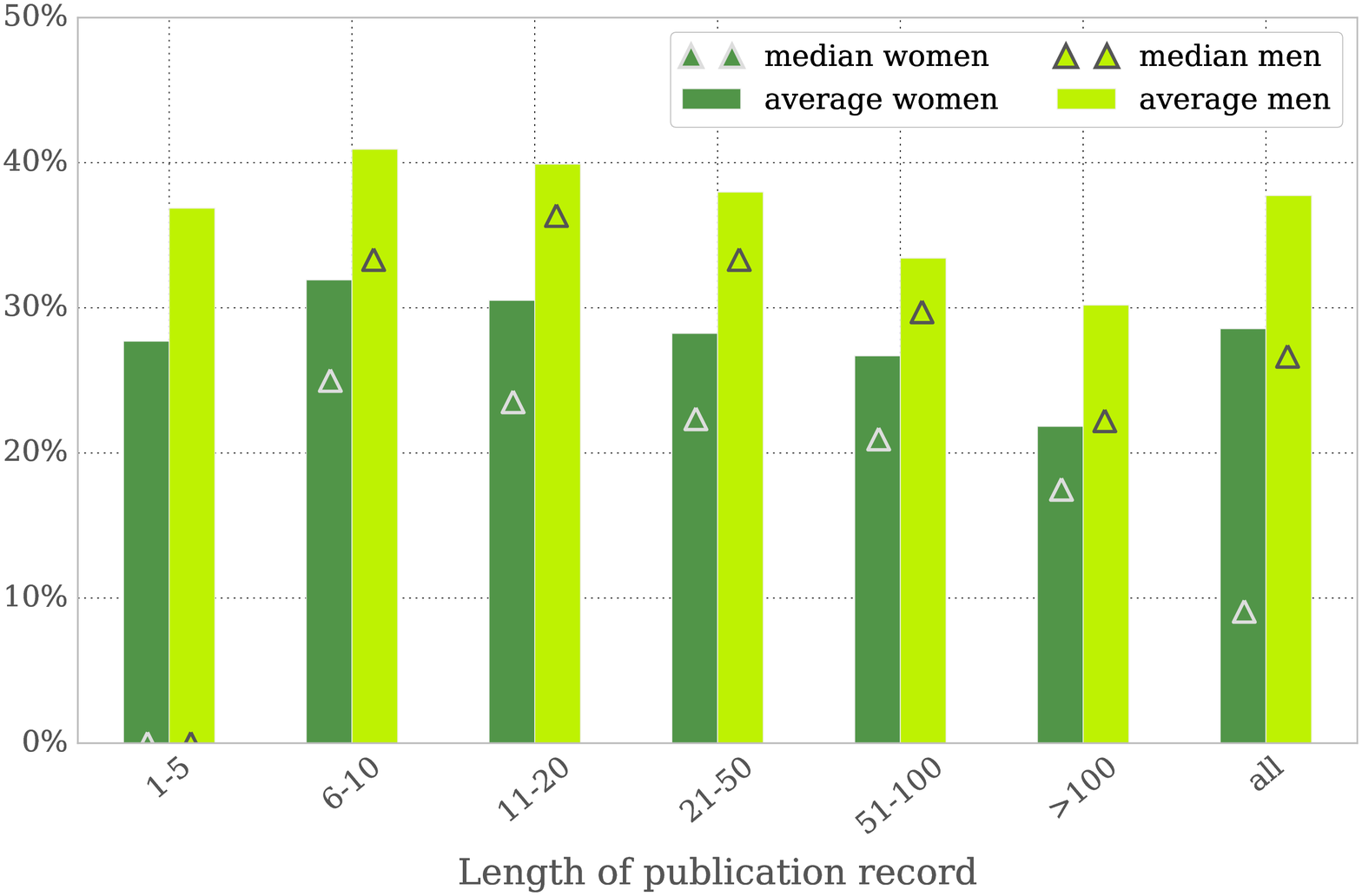}
		\caption{{\bf Mean and median of the amount of single-authored publications per author.} 
		Authors are grouped in segments defined by length of publication record: at most 5, 10, 20, 50, or more than 100 publications. The rightmost bar refers to the full data set.}
		\label{Fig 8}
	\end{center}
\end{figure}

On average, men write 38\% of their scientific records as single authors, in contrast to 29\% among women. This trend remains unchanged when looking at groups according to the total number of publications, in segments of at most 5, 10, 20, 50, or more than 100. Fig~\ref{Fig 8} illustrates this conclusion, showing a statistically significant gap between men and women in the percentage of single-authored publications across all segments. A similar tendency has also been reported in~\cite{Schucan-Bird.K:2011} for UK-authored journal articles in Social Sciences, or in~\cite{Cunningham.SJ.Dillon.SM:1997} for the field of Information Systems; in the later study, the percentage of male authors with a single-authored publication was twice as large as for their female counterparts. 

Sooner or later, many mathematicians publish an article on their own. In our data set, this is the case for 52\% of the women; in the men's sample this holds true for 61\%. The single-authored paper is usually not the author's first publication, in particular in the case of female mathematicians. Women write their first publication most often together with coauthors: only 33\% of their first publications are single authored, for men it is 43\%. 

\paragraph*{Women's and men's coauthor networks are of similar size}

\begin{figure}[h!]
	\begin{center}
		\includegraphics[width=\columnwidth]{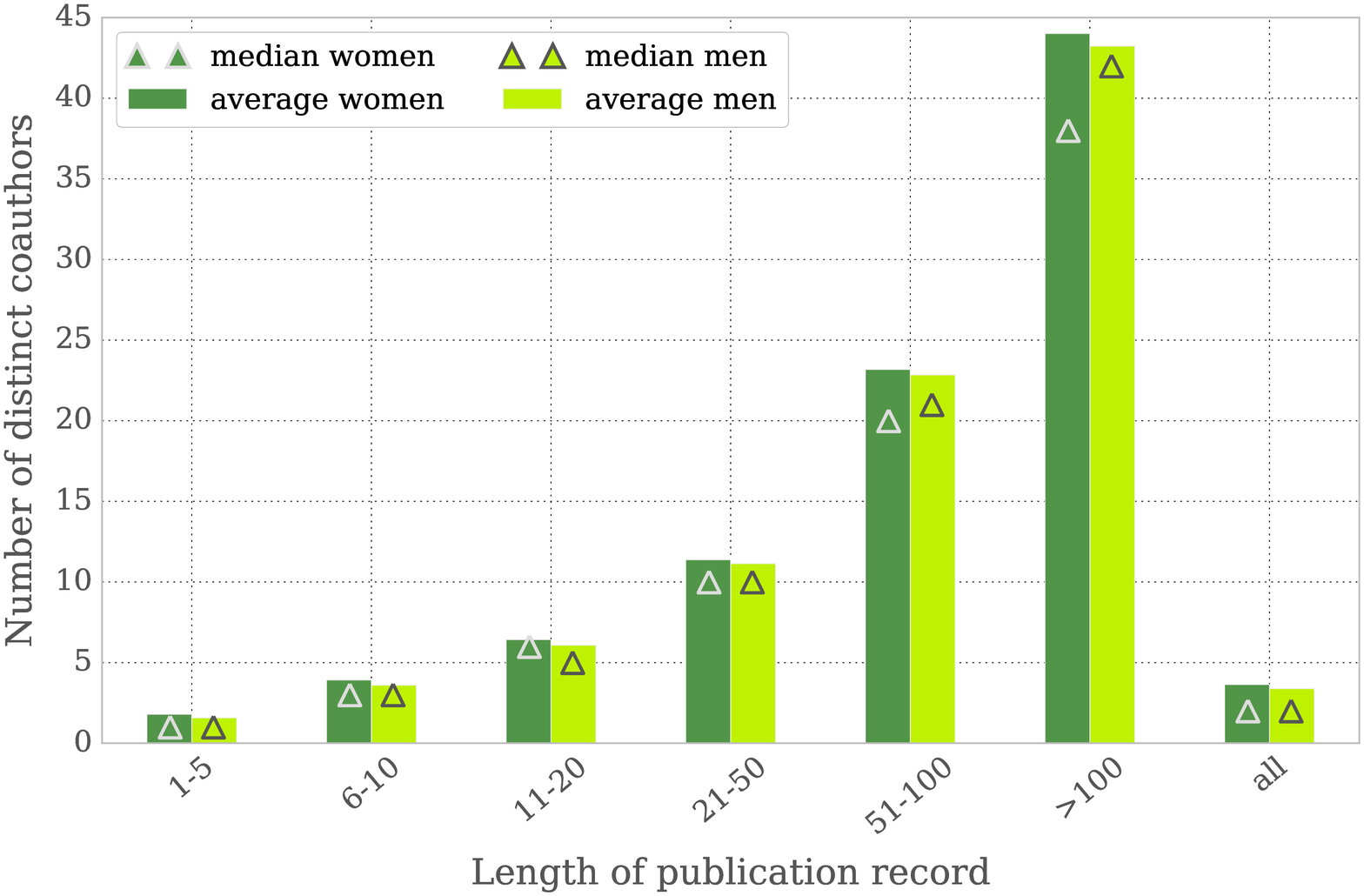}
		\caption{{\bf Mean and median of the number of distinct coauthors per author.} 
		Authors are grouped as in Fig~\ref{Fig 8}.}
		\label{Fig 9}
	\end{center}
\end{figure}

Although women mathematicians have significantly less single-authored publications than their male counterparts, they do not collaborate with a considerably larger number of peers. On average, women have 3.6 distinct coauthors, 0.2 more than men.
Fig~\ref{Fig 9} shows the network size for men and women in more detail, split according to the length of publication record. The dark and light green bars display the average network sizes of women and men, respectively, with the enclosed triangles representing the medians. The mean network size of women is slightly larger in each of the groups; however, the gap decreases with a growing number of publications from factor 1.15 in the group of authors with at most 5 publications (i.e. the network size of women in this group is on average 1.15 times larger than that of men) to 1.02 in the two groups of authors with more than 50 publications. The median values for women and men are equal in the first two and the fourth group; for authors with 11-20 publications in total, women have a median of 6 compared to a median of 5 for men. In the last two  groups, the median value for men is larger; among authors with more than 100 publications the group of men exhibits a median of 42 versus 38 for women. 

Both measures suggest that women and men collaborate with a similar number of peers. However, when measuring the network size with respect to the time line instead of the total number of publications, then Fig~\ref{Fig 3} suggests that the networks of women are even smaller than those of men, at least when looking at an early career stage, since women publish around 10\% less articles than men within the first 5 years of their research~\cite{Pierson:2014}.

Bozeman and Gaughan~\cite{Bozeman.B.Gaughan.M:2011} statistically analysed collaboration patterns based on surveys of STEM researchers in the United States. In contrast to the authors' expectations, women had a slightly (but not significantly) larger collaboration network. However, the study did not include information on publications, and hence could not elaborate on the relation between the number of collaborators and single-authored publications.

It is a common assumption that a larger coauthor network is indicative of the visibility of an author's research. For instance, Bentley~\cite{Bentley.P:2012} shows for Australian university research that international collaborators were one of the main factors associated with increased publication productivity among both men and women. Our findings support the conclusion that men and women collaborate with similarly-sized peer groups. However, women seem to do so at the expense of their own single-authored publications, which might undermine their profiles as fully independent researchers.

\subsection*{Topical distribution}

We study the distribution of authors and publications in mathematical fields as represented by the MSC 2010 (see subsection \textit{The zbMATH corpus}). We focus on the first and second MSC levels; the first level represents mathematical fields such as ``Algebraic Topology'' or ``Dynamical Systems'', whereas the second level allows for a more granular view of related mathematical communitites. We restricted to the primary MSC code assigned to an article since this usually captures the main content of the research.

\paragraph*{Women's distribution across research fields is not homogeneous}

Women publish in all of the 63 mathematical fields represented by the first hierarchy level of the MSC; out of the 528 topic-based subfields at the second MSC-level, there are only three in which none of the women in our data set have published either alone or with a coauthor. By restricting our search to single-authored publications, we find 13 second-level classes in which women have not yet published, such as ``$K$-theory in number theory'' (19F) or ``Function theory on the disk'' (30J).

\begin{figure}[h]
	\begin{center}
		\includegraphics[width=\columnwidth]{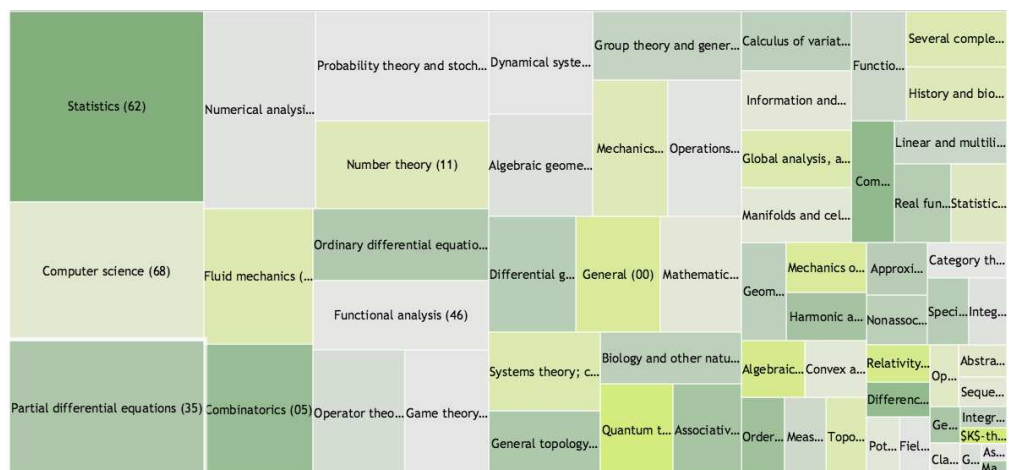}
		\caption{{\bf Distribution of women authors across mathematical fields and subfields.}
		Mathematical subjects are represented by MSC classes as a hierarchical map. Data includes publications of all mathematicians classified as women or men. An intense green color stands for a relatively high share of women's publications, an intense yellow for a relatively low share. The size of the rectangle is proportional to the field size. Under \url{https://doi.org/10.5281/zenodo.51147}, the detailed numbers are shown when hovering over the rectangles. Clicking on a rectangle shows the distribution on the second MSC-level in that field; a right-click returns to the first-level view.}
		\label{Fig 10}
	\end{center}
\end{figure}

However, women are not equally represented across mathematical areas. Using hierarchical maps as displayed in Fig~\ref{Fig 10} we visualize the share of women resp. women's publications across the MSC classes on level 1 and 2: each rectangle represents a mathematical class, its size is proportional to the number of authors resp. publications in that field, and the darker the color, the larger their share. 

This confirms the results by West~et al.~\cite{West.JD.Jacquet.J.ea:2013, West.JD.Bergstrom.CT.ea:2012} who showed an uneven distribution based on publications from JSTOR. Unfortunately, we cannot compare the results in detail since West~et al. did not use a classification schema such as the MSC but instead employed an unsupervised hierarchical clustering approach to assign publications to topics. 

\begin{figure}[h]
	\begin{center}
		\includegraphics[width=\columnwidth]{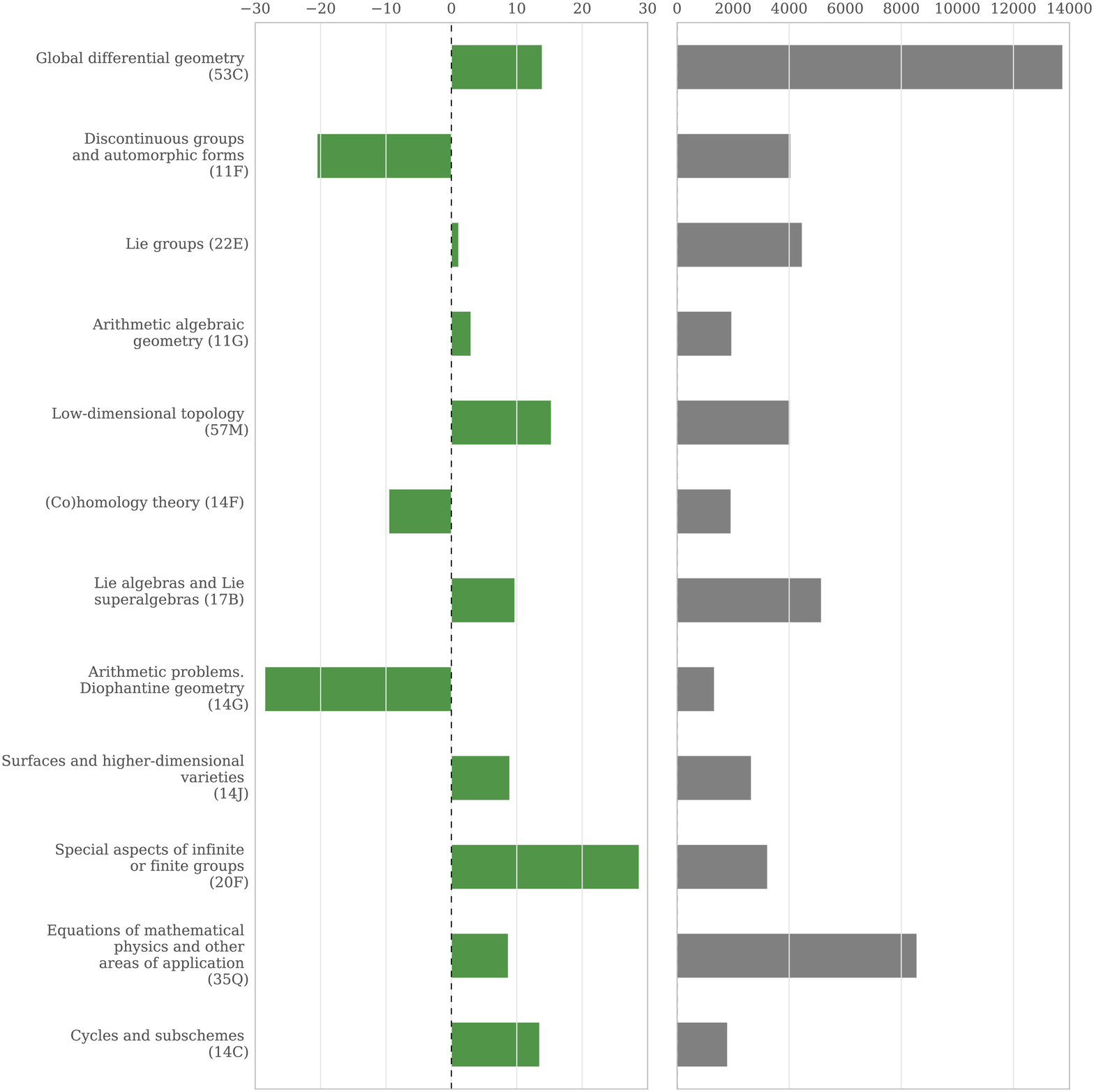}
		\caption{{\bf Single-authored publications in twelve research fields.} 
		These second MSC-level subjects dominate the articles in the analyzed three top journals. The bar chart on the left shows the deviation of the share of women's publications from the average, on the right-hand side is the size of the research field in terms of the overall number of single-authored publications.}
		\label{Fig 11}
	\end{center}
\end{figure}

In subsection \textit{Perception of quality} we showed that the presence of women in three of the most renowned mathematical journals---the ``Annals of Mathematics'', ``Inventiones Mathematicae'' and the ``Journal of the American Mathematical Society''---is low, and reveals no significant improvement over the last three decades. The three journals, despite aiming to cover all fields in mathematics, display a bias towards certain fields. For instance, the five areas ``Algebraic Geometry (14)'', ``Number Theory (11)'', ``Differential Geometry (53)'', ``Dynamical Systems and Ergodic Theory (37)'' and ``Group Theory and Generalizations (20)'', already cover more than 50\% of all publications in these journals. On the second MSC level we see a similar dominance of certain subfields, leading to a coverage of more than a third of all publications by only 12 classes which are shown in Fig~\ref{Fig 11}. To investigate whether the selection of research areas might be responsible for the underrepresentation of women in these journals, we checked for correlation between a particular field's presence in the three journals and the women's share in this field. It turns out that women are well represented in the top 12 level-2 classes ($\sim\!9\%$ above the average when considering all publications and $\sim\!13\%$ when restricting to single-authored publications). When looking at publications in those fields instead of authors, the numbers even out ($\sim\!-3\%$ when looking at all publications and $\sim\!4\%$ when restricting to single-authored publications). Most level-1 fields with strong dominance in the three journals show an average or even higher presence of women, with the traditional field of ``Number Theory'' being the most significant negative exception. 

\paragraph*{Women's distribution of research topics is narrower than men's}

\begin{figure}[h]
	\begin{center}
		\includegraphics[width=\columnwidth]{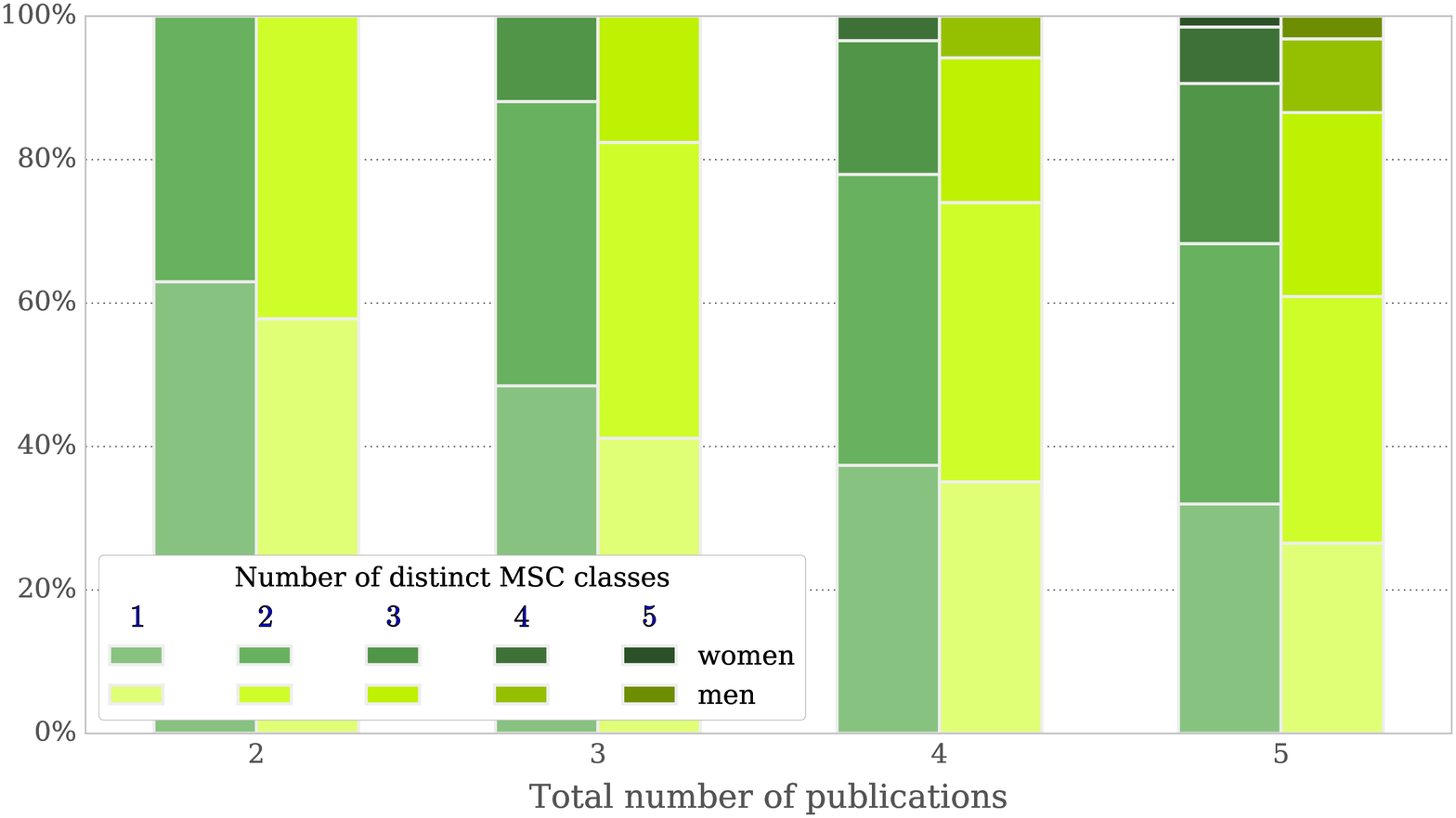}
		\caption{{\bf Diversity of research fields for women and men with 2 to 5 publications in total.}
		For a fixed total number of publications ranging from 2 to 5, the same number of female and male authors has been sampled and the number of different MSC classes per author has been computed.}
		\label{Fig 12}
	\end{center}
\end{figure}

As shown in subsection \textit{Perception of quality}, women and men collaborate with a comparable number of peers when measuring against groups of authors with a similarly-sized publication record. Taking the point in time as scale instead of the total number of publications (see Fig~\ref{Fig 3}), the networks of women are smaller than those of men. We see a similar pattern when considering the variation of topics: as shown in Fig~\ref{Fig 12}, publications of women spread across less MSC classes than those of men, while keeping the overall number of publications per author fixed. Their stronger focus on a smaller set of research topics could potentially be another drawback for women in academia, since it reduces the already low number of possible faculty positions suitable for applying to. 

\section*{Discussion}

Our study is based on comprehensive data from publications in all areas of mathematics research between 1970 and 2013. Backed by a statistically significant data corpus and scanning a rather broad set of determinants for scholarly productivity, we audit the status of women mathematicians in the publication landscape and identify research questions which would provide additional substantial insights. 

We have split all authors in our data set in cohorts according to the year of their first publication and have investigated their progression 5 and 10 years into their careers. We have looked into variations in publication patterns between men and women regarding presumed quality of the journals they publish in, their coauthorship networks, and the topical distribution of their scholarly output. We report statistically significant differences among genders in all the aforestated aspects. 

The share of women mathematicians who have been active for at least 5 to 10 years has roughly tripled since the 1970s (see Fig~\ref{Fig 2}). However, women abandon their academic careers within a decade of having started them at a much higher rate than men (see Fig~\ref{Fig 4}). This indicates the existence of a ``leaky pipeline'' that slows down the inclusion of women in the higher ranks of academia, and it invalidates the hypothesis that, given the steady increase of women entering mathematics at the undergraduate level, their incorporation to higher academic ranks should simply follow suit. Furthermore, the so-called productivity gap, stating that women produce less publications within a fixed period of time, is clearly visible in every cohort. However, it has been narrowing down drastically since 1970 with a visibly positive trend in recent years.

The perceived quality of the journal in which an article is published is one of the most common means of evaluating a researcher's academic output. There is legitimate and compelling criticism of this approach. Beyond all objections to particular shortcomings of these measurements, what they all have in common is the impossibility of providing any information on individual articles. While we believe that journal rankings are no proper means to evaluate scientists and their work, they are still widely adopted to measure research quality. With respect to such a measurement scheme,  women are extremely underrepresented in top-ranked journals. A negative correlation exists between journal rank, measured using the ERA ranking or the journal impact factor, and the percentage of women among its authors, a pattern that remains fairly stable over time. A further investigation of three highly renowned mathematics journals shows an unusually low percentage of published women, with almost no improvement over time. Several explanations may be proposed for this fact, from biased editorial boards and reviewers to lack of strategic publication choices on the part of women. Using our data, we were able to exclude the choice of the research field as possible reason, since women are well represented in the fields favoured by the considered journals. Certainly, it would be worth to further investigate if and why women simply submit less articles to prominent journals or whether their papers are more often rejected. Nonetheless, our findings uncover a specific phenomenon correlating journal rank and scarcity of published women, and the reasons for these differences need to be studied in more detail.

With respect to collaboration through coauthorship, our results indicate that men publish far more single-authored papers than women. Since a solid record of single-authored publications seemingly shows one's ability of independence and creativity and is hence associated with ``scientific autonomy",  the disparity could have a negative effect on women's career perspectives, in particular at an early stage. Considering that their network is not significantly larger, this suggests that women regularly collaborate with the same researchers and begs the question of \textit{with whom} women and men collaborate. Are the most frequent coauthors the doctoral advisor and close collaborators? Likewise, an analysis of the coauthors' gender could help understand whether gender induces distinctive publication dynamics, e.g in an all-female versus a mixed-gender collaboration. Additionally, it would be interesting to explore if factors such as the academic position or the country/internationality of the collaborators have an effect on the collaboration patterns.

It seems sensible to consider the aforementioned questions also in a field-focused setting since the presence of women strongly deviates across research area, as shown in Section \textit{Topical distribution}. It would also be interesting to explore in more detail which factors motivate women to enter a certain research area. What constitutes an inviting ``atmosphere" within a field, what is the importance of role models, mentoring, etc.? How does this relate with the age and interdisciplinarity of a field? Do new and emerging areas produce a different mind-set and are hence more attractive for women scientists?   

We note that our analysis is robust with respect to the counting scheme. In certain cases, e.g. when papers authored by men have a larger number of authors than papers by women, it might be more suitable to analyse distributions of fractional authorships where an authorship is weighted reciprocally to the overall number of authors of the paper~\cite{Lariviere2013,Paul-Hus.A.Bouvier.RL.ea:2015}. For our data set, however, both full and fractional counting lead to very similar results. This is explained by the dominance of single-authored publications in mathematics and the relative balance of women's and men's network sizes. Since women in our data set have, on average, more coauthors than men, the gender gap is even slightly stronger in terms of fractional authorships. 

Our findings support the thesis that gender-related publication patterns exist and are one of the factors that lead to an underrepresentation of women in mathematics. Our analysis cannot give a definite answer to the main motivating question, i.e. \textit{why} the gender gap in mathematics. The differences in research productivity according to gender are by no means trivial to measure or explain. Likewise, research productivity alone does not fully translate into career success, however defined. Various studies, including those cited in Section \textit{Related work}, have looked at the factors that might contribute to a higher or lower research output, including demographics, family-related circumstances, working environment, etc.  
Furthermore, academic excellence is an inherently gendered concept~\cite{Brink.M.Benschop.Y:2012}, with several aspects contributing to its fabrication. Even if the publication and collaboration patterns of men and women showed no significant differences, a completely gender-blind approach to promotion in mathematics academic positions is unfortunately not to be expected. Hence a thorough discussion on the topic must necessarily expand towards further aspects. Perhaps a first one is the question of who or what is setting the pace and norms of scholarly publishing. In the work of Brink and Benschop~\cite{Brink.M.Benschop.Y:2012} the concept of ``masculine career path'' (p. 513) is put forward: every researcher (man or woman) whose publishing record deviates from the traditional trajectory, maybe due to career interruptions or periods of part-time work, is at dire disadvantage for promotion. A similar remark is found in Bentley~\cite{Bentley.P:2012}, which mentions the institutionalization of male dominance by associating merit primarily with scholarly output, an activity for which men receive greater patronage and support. The lack of encouragement that many female doctoral students perceive when it comes to applying for funding or to collaborate on research projects already sets back their careers. Women are forced to either ``catch up'' later or, or are driven out of the field altogether at a higher rate than men, as evidenced by our cohort analysis. The fact that in many cases the number of publications is the main factor in defining the quality of research constitutes a drawback for the numerous women with nonlinear career trajectories or less research time.

Using one of the largest metadata sources on mathematical publications we have shown that the collective actions of all actors in mathematics publishing result in significant differences between two gender groups regarding key factors of academic careers.
A more thorough and detailed analysis of single aspects would clearly help to explain the observed disparities between women and men. For this, additional data such as affiliations, acceptance and rejection statistics from journals, or the country of origin would be very helpful; the latter could also be used for an improvement of the gender assignment algorithm which constitutes the base of our study. In fact, we would prefer not to use such an assignment schema at all but ask the authors instead how they identify.  
Our analyses have unveiled significant differences which may already foster individual takes on publishing, hence career decisions, and inspire institutional measure. We believe that the gained insights will deliver input for further discussions within the community for potential collective and/or individual strategies to overcome the current situation.

\section*{Acknowledgements}
We gratefully acknowledge the support of FIZ Karlsruhe for granting us access to their publication data. We thank the European Women in Mathematics network and specially its former convenor Marie-Françoise Roy for feedback in the first stages of this project and useful comments on an earlier paper draft.

\end{document}